\documentclass[a4paper, 11pt]{amsart}

\usepackage{enumerate}
\usepackage{mathrsfs}
\usepackage{amssymb}
\usepackage[all]{xy}

\title[Okounkov bodies for ample line bundles]
{Okounkov bodies for ample line bundles}
\author{Henrik Sepp\"anen*} 

\thanks{*supported by the DFG Priority Programme 1388 ``Representation 
Theory"}

\date{\today}

\address{Henrik Sepp\"{a}nen,
Mathematisches Institut,
Georg-August-Universit\"at G\"ottingen,
Bunsenstra\ss e 3-5, 
D-37073 G\"ottingen,
Germany}
\email{hseppaen@uni-math.gwdg.de}

\newcommand{\R}{\mathbb{R}}
\newcommand{\N}{\mathbb{N}}

\newtheorem{prop}{Proposition}[section]
\newtheorem{lemma}[prop]{Lemma}

\newtheorem{thm}[prop]{Theorem}
\theoremstyle{definition}

\newtheorem{rem}[prop]{Remark}

\begin{document}
 
\maketitle

\begin{abstract}
Let $\mathscr{L} \rightarrow X$ be an ample line 
bundle over a nonsingular complex projective variety $X$. 
We construct an admissable flag 
$X_0 \subseteq X_1 \subseteq \cdots \subseteq X_n=X$ of 
subvarieties for which the associated Okounkov body for 
$\mathscr{L}$ is a rational polytope.

\end{abstract}

\section{Introduction}
The \emph{Okounkov body} of a line bundle 
$\mathscr{L} \rightarrow X$ over a nonsingular $n$-dimensional complex 
variety is a convex set in $\R^n$ which carries information about the
section ring $\bigoplus_{k=0}^\infty H^0(X,\mathscr{L}^k)$ of $\mathscr{L}$. 
The idea is to construct a valuation-like function (meaning that 
it has the properties of a valuation in the ring-theoretic sense,
even if it only defined on homogeneous elements)
$$v: \bigsqcup_{k \in \N} H^0(X,\mathscr{L}^k)\setminus \{0\}  \rightarrow \N_0^n,$$ 
and define the semigroup 
\begin{equation*}
S(\mathscr{L},v):=\{(k,v(s)) \mid s \in H^0(X,\mathscr{L}^k) \setminus \{0\}\}. 
\end{equation*}
The associated Okounkov body is then defined as 
\begin{equation}
\Delta(\mathscr{L},v):=
\overline{\mbox{conv}}\left\{\frac{1}{k} (k,v(s)) \mid  s \in H^0(X,\mathscr{L}^k) 
\setminus \{0\}\right\}.
\end{equation}
Roughly speaking, the function $v$ will be defined as the successive orders 
of vanishing of sections along a flag 
$X_0 \subseteq X_{1} \subseteq \cdots \subseteq X_{n-1} \subseteq X_n=X$ 
of irreducible nonsingular subvarieties such that $X_i$ has dimension $i$.

These bodies were introduced by Okounkov in \cite{Ok96} from a 
representation-theoretic point of view. However, he only considered the 
semigroup defined by the values of $U$-invariants, where $U$ is the unipotent
radical of a Borel subgroup, $B$, of $G$. It turned out that
the elements of height $k$ in this semigroup carry information about 
the decomposition of the $k$th piece, and that the Euclidean volumes 
certain slices of the associated convex body describe the
asymptotics of the decomposition.

It has since then become an interesting problem \emph{per se} to study the 
semigroups and associated convex bodies defined by considering the values 
of all sections of a line bundle, without the assumption of
a group action (\cite{LM09}), and their relation to the 
geometry of the line bundle $\mathscr{L}$.
One crucial connection between the section ring of $\mathscr{L}$ and the
Okounkov body $\Delta(\mathscr{L},v)$ is the identity
\begin{equation}
\mbox{vol}\Delta(\mathscr{L},v)=\lim_{k \rightarrow \infty}
\frac{\mbox{dim} H^0(X,\mathscr{L}^k)}{k^n},
\end{equation}
which holds for the linear series of a big line bundle over
a smooth irreducible projective variety $X$ (see \cite{LM09}).

A fundamental question to ask is whether $\Delta(\mathscr{L},v)$ is 
a rational polytope, i.e., the convex hull of finitely many 
points in $\mathbb{Q}^n$. This was conjectured by Okounkov
in \cite{Ok96} in his setting where he only considered $U$-invariants. 
It is true in some cases, e.g., equivariant line 
bundles over toric varieties. This follows from recent results 
by Kaveh and Khovanskii where they consider a setting that
generalizes Okounkov's, namely group representations on 
graded algebras of meromorphic functions on $X$, cf. \cite{KK10}.
In fact, the Okounkov body of a torus-equivariant line
bundle over a projective toric variety equals the \emph{moment polytope} 
(cf. \cite{Br86}) associated to the graded representation.
Recently Kaveh (\cite{Kav11}) has studied this problem for a 
line bundle over a full flag variety $G/B$. He considers a 
flag of translated Schubert varieties and shows, 
using Bott-Samelson resolutions, that the corresponding 
Okounkov body is a polytope by identifying it with a 
Littelmann string polyope.

It is however not true in general that $\Delta(L,v)$ is 
a rational polytope. A counterexample can be found in 
\cite[Section 6.3]{LM09}.

In this short note, we prove that if $\mathscr{L}$ is ample, there
exists a flag of nonsingular irreducible subvarieties
$X_0 \subseteq X_1 \subseteq \cdots \subseteq X_n:=X$ for which the 
associated Okounkov body is a rational polytope.

We would like to point out that the result in this note
also follows from a more general result obtained by
Anderson, K\"uronya and Lozovanu in \cite{akl}, where they study 
the problem of finding a finitely generated semigroup defining the 
Okounkov body.

\noindent {\bf Acknowledgement}. 
We would like to thank Dave Anderson, Alex K\"uronya and 
Robert Lazarsfeld for helpful discussions on this topic.

\section{A polyhedral Okounkov body}
Let $X$ be a nonsingular $n$-dimensional complex variety and 
let $\mathscr{L} \rightarrow X$ be an ample line bundle. 
Since Okounkov bodies satisfy the property $\Delta(\mathscr{L}^k)=k\Delta(\mathscr{L})$, 
$k \in \N$ (cf. \cite[Proposition 4.1.]{LM09}), we may without loss of
generality assume that $\mathscr{L}$ is very ample.

By Bertini's theorem there exist sections 
$\xi_1,\ldots, \xi_{n-1} \in H^0(X, \mathscr{L})$ such that 
the zero sets
\begin{align}
X_i:=\{x \in X \mid \xi_1(x)=\cdots =\xi_{n-i}(x)=0\}, \quad i=1,\ldots, n-1
\end{align}
define nonsingular, irreducible subvarieties of $X$. Let 
$X_0:=\{p\} \subseteq X_1$ for some point $p \in X_1$.
Let 
\begin{align*}
v: \bigsqcup_{k \geq 0} H^0(X, \mathscr{L}^k) \setminus \{0\}
\rightarrow \N^n
\end{align*}
be the valuation-like function defined by 
the flag 
\begin{align}
X_0 \subseteq X_1\subseteq \cdots \subseteq X_{n-1} \subseteq X_n:=X,
\label{E: flag}
\end{align}
and let $\Delta \subseteq \R^n$ be the 
associated Okounkov body. For the construction, we refer to \cite{LM09}.
Similarly, let 
\begin{align*}
v_1: \bigsqcup_{k \geq 0} H^0(X_1, \mathscr{L}^k\mid_{X_1}) \setminus \{0\}
\rightarrow \N
\end{align*}
be the valuation-like function defined by the flag $X_0 \subseteq X_1$,
and let $\Delta_1 \subseteq \R$ be the corresponding Okounkov body.
Then 
\begin{align*}
\Delta_1=[0, b] \subseteq \R,
\end{align*} 
where $b=\mbox{deg} \,\mathscr{L}\mid_{X_1}$ (cf. \cite[Ex. 1.14]{LM09}).
\begin{lemma} \label{L: okincl}
The inclusion 
\begin{align*}
[0,b] \times \{0\}\times \cdots \times \{0\} \subseteq \Delta
\end{align*}
holds.
\end{lemma}

\begin{proof}
By the closedness of $\Delta$ it suffices to prove that 
$(c,0,\ldots, 0) \in \Delta$ for every $0<c<b$. 

By the ampleness of 
$\mathscr{L}$ there exists an $N_0 \in \N$ such that 
the restriction maps
\begin{align*}
R_j: H^0(X, \mathscr{L}^j) \rightarrow H^0(X_1, \mathscr{L}^j\mid_{X_1}) 
\end{align*}
are surjective for $j \geq N_0$. 
If now $c \in (0,b)$, choose a point $v_1(\tau)/m \in (c, b)$ for a section 
$\tau \in H^0(X_1, \mathscr{L}^m\mid_{X_1})$.
Hence, for sufficiently large $N \in \N$,
there exists a section
$\widetilde{\tau_N} \in H^0(X, \mathscr{L}^{Nm})$
with $R_{Nm}(\widetilde{\tau_N})=\tau^N$. 
Then $v(\widetilde{\tau_N})=(v_1(\tau^N),0,\ldots, 0)=(Nv_1(\tau),0,\ldots, 0)$, 
so that $(v_1(\tau)/m,0,\ldots, 0) \in \Delta$. By the convexity of 
$\Delta$ we then have $(c,0,\ldots, 0) \in \Delta$.
\end{proof}

\begin{rem}
The lemma would of course also hold if the subvarieties $X_i$ 
were defined by unequal line bundles $\mathscr{L}_i$, $i=1,\ldots, n-1$. 
For a related result for line bundles $\mathscr{L}$ that are not ample, 
but merely big, cf. \cite[Theorem B]{jow}.
\end{rem}

Let $e_1,\ldots, e_n$ be the standard basis for $\R^n$.
\begin{thm} \label{T: okamp}
The Okounkov body $\Delta$ is the
convex hull of the set 
$$\{0, be_1, e_2, \ldots, e_n\}.$$
\end{thm}

\begin{proof}
First of all, $v(\xi_{n-i+1})=e_i$ for $i=2,\ldots, n$. 
Since the linear system defined by $\mathscr{L}$ is basepoint-free, 
$0 \in v(H^0(X, \mathscr{L}) \setminus \{0\})$. By Lemma \ref{L: okincl}
$be_1 \in \Delta$. Hence the convex hull of the points
$0,be_1, e_2,\ldots, e_n$ is a subset of $\Delta$.

It thus suffices to prove that for every $a=v(s)$ for some 
$s \in H^0(X, \mathscr{L}^k) \setminus \{0\}$ there exist 
$x_0,\ldots, x_{n} \geq 0$ 
such that 
\begin{align}
&x_0+\cdots+x_{n}=k, \nonumber\\ 
&a=x_0\cdot 0+x_1 be_1+x_2e_2+\cdots +x_{n}e_n. 
\label{E: convcomb}
\end{align}
For this, let $a=(a_1,\ldots, a_m,0,\ldots, 0)=v(s)$ for 
$s \in H^0(X, \mathscr{L}^k)$ and assume that $a_m\geq 1$. 
Then $(s/\xi_{n-m+1}^{a_m})\mid_{X_{m-1}}$ defines a section in \\
$H^0(X_{m-1}, \mathscr{L}^{k-a_m}\mid_{X_{m-1}})$, and by iteration 
we can write 
\begin{align*}
a=(a_1,0,\ldots, 0)+\sum_{i=2}^m a_iv(\xi_{n-i+1}),
\end{align*} 
where 
$a_1=v_1(t)$ for the section 
$$t=(s/(\xi_{n-1}^{a_2}\cdots \xi_{n-m+1}^{a_m}))\mid_{X_1} \in 
H^0(X_1, \mathscr{L}^{k-\sum_{i=2}^ma_i}\mid_{X_1}).$$
Since the line bundle $\mathscr{L}^{k-\sum_{i=2}^ma_i}\mid_{X_1}$ is 
effective we must have $p:=k-\sum_{i=2}^ma_i \geq 0$. Indeed, 
$\mathscr{L}^r\mid_{X_1}$ is effective for every $r \geq 0$, 
and hence $\mathscr{L}^r \mid_{X_1}$ cannot be effective for 
any $r<0$. 
Since $\Delta_1=[0,b]$ there exist real $x_0, x_1 \geq 0$ 
with $x_0+x_1=p$ and 
\begin{align}
a_1=x_0 \cdot 0+x_1 b. \label{E: conv1},
\end{align}
and hence
\begin{align*}
a_1e_1=x_0 \cdot 0+x_1be_1.
\end{align*}
By putting 
\begin{align*}
x_i:=a_i \quad i=2,\ldots, n,
\end{align*}
we have thus found nonnegative real numbers $x_0,\ldots, x_n$ satisfying 
\eqref{E: convcomb}.
\end{proof}

\section{Concluding remarks}
The subvarieties $X_i$ occurring in the flag \eqref{E: flag} are unfortunately only 
defined very implicitly. It would be interesting to find more explicit examples of 
admissable flags that lead to polyhedral Okounkov bodies. In particular, in the original 
setting developed by Okounkov, involving the presence of a group action, it is desirable 
to have an admissable flag with a certain invariance property in order to relate the 
Okounkov body to asymptotics of multiplicities for irreducible representations.

Another interesting problem is to find admissable flags yielding polyhedral 
Okounkov bodies for a larger class of line bundles, or even polyhedral global 
Okounkov bodies (cf. \cite[Section 4.2]{LM09}).

\end{document}